\documentclass[10pt]{amsart}

\allowdisplaybreaks
\usepackage{amsmath}
\usepackage{amsthm}
\usepackage{amssymb}
\usepackage{amsfonts}
\usepackage{latexsym}
\usepackage{array}
\usepackage{verbatim}
\usepackage{mathptmx} 
\usepackage{graphicx} 
\usepackage{color}
\usepackage{mysects}  
\usepackage{fancyhdr} 
\usepackage{enumerate}
\usepackage[labelformat=simple]{caption}
\usepackage{setspace}
\usepackage{microtype}
\usepackage{bm}

\title{The number of roots of polynomials of large degree in a prime field.}
\author{Amit Ghosh and Kenneth Ward}
\newtheorem*{theorem1}{Theorem 1}
\newtheorem*{theorem2}{Theorem 2}
\newtheorem*{theorem3}{Theorem 3}
\newtheorem*{lemma1}{Lemma 1}
\newtheorem*{lemma2}{Lemma 2}
\newtheorem*{lemma3}{Lemma 3}
\newtheorem*{lemma4}{Lemma 4}
\newtheorem*{lemma5}{Lemma 5}
\newtheorem*{lemma6}{Lemma 6}
\newtheorem*{lemma7}{Lemma 7}
\newtheorem*{rem}{Remark}
\theoremstyle{definition}
\theoremstyle{definition}
\newtheorem*{definition}{Definition}

\begin{document}
\dedicatory{Dedicated to Roger Heath-Brown on the occasion of his sixtieth birthday.}
\begin{abstract} We establish asymptotic upper bounds on the number of zeros modulo $p$  of certain polynomials with integer coefficients, with $p$ prime numbers arbitrarily large.  The polynomials we consider have degree of size $p$ and are obtained by truncating certain power series with rational coefficients that satisfy simple differential equations.  \end{abstract}

\maketitle
\thispagestyle{empty}
\pagestyle{fancy}
\fancyhead{}
\fancyhead[LE]{\thepage}
\fancyhead[RO]{\thepage}
\fancyhead[CO]{\small The number of roots of polynomials of large degree in a prime field}
\fancyhead[CE]{\small Amit Ghosh and Kenneth Ward}
\renewcommand{\headrulewidth}{0pt}
\headsep = 1.0cm
\fancyfoot[C]{}
\vskip 1cm
\renewcommand{\thefootnote}{ } 
\renewcommand{\footnoterule}{{\hrule}\vspace{3.5pt}} 


\section{Introduction}

In 1996, Heath-Brown \cite{He1} showed that for sufficiently large prime numbers $p$,  the Heilbronn sum $\sum_{(x\mod p)} \exp(2\pi i \frac{a x^p}{p^2}) $ is bounded above by $Cp^\frac{11}{12}$ when  $a$ is not divisible by $p$, so that Weyl's criterion implies the uniform distribution of the sequence $\{x^p \mod p^2: x \mod p\}$ (here $C$ denotes an absolute constant). The novelty of his approach was to reduce the problem (after an application of Cauchy-Schwarz and a suitable change of variables) to counting the number of solutions to a polynomial congruence modulo $p$ and then to achieve a non-trivial count by a modification of Stepanov's method appearing in the proof of the Riemann hypothesis for a curve over a finite field \cite{St, Bo}. Coincidentally, in 1992, Mit'kin \cite{Mi} had considered just this counting question for two polynomials (one of which happened to be the one used in \cite{He1}). Moreover, he used the same methods and obtained the same bound of \cite{He1}. Our focus is to consider generalisations of this counting problem.

 For a prime number $p\geq 3$,  we consider polynomials $F(x)$ having rational coefficients with denominators coprime to $p$ and with degree $d$.  Let $N_{p}(F)$ denote the number of  solutions to the congruence $F(x) \equiv 0 \mod p$. We are interested in bounds for $N_{p}(F)$ with $p$ sufficiently large when the degree $d$ has size proportional to $p$. It is clear that in this generality one cannot say better than the trivial bound $N_{p}(F) \ll p$. Indeed if $F(x) = x^p -x$, then $N_p(F) = p$ and moreover for any $a \not\equiv 0 \mod p$, $N_p(F-a)=0$. A much deeper example can be found in the theory of elliptic curves. Consider the Legendre elliptic curve over ${\mathbb F}_p$ for $p>3$, given by $y^2 = x(x-1)(x-\lambda)$ with $\lambda \neq 0,1$. Let $H_p(\lambda) = \sum_{n=0}^{D} \binom{D}{n}^{2}\lambda^n$ with $D=\frac{p-1}{2}$, be the Hasse-Deuring invariant. Finally let $a_p(\lambda) = N-p+1$ where $N$ counts the number of ${\mathbb F}_p$-rational points on the curve. Then it is known (Igusa \cite{Ig58}, Dwork \cite{Dw62}) that $a_p(\lambda) \equiv (-1)^D H_p(\lambda) \mod p$. It therefore follows that the number of solutions to the polynomial congruence $H_p(\lambda) \equiv A \mod p$, for any fixed A is essentially the same as the number of isomorphism classes of such elliptic curves satisfying the condition $a_p(\lambda)\equiv (-1)^D A \mod p$. The Hasse estimates imply that there are no solutions if $|A|> 2\sqrt{p}$. On the other hand, if $|A| < 2\sqrt{p}$, Deuring showed that the number of such classes is essentially the Kronecker class number $H(\lambda^2 -4p)$ of an imaginary quadratic field. It can then be shown that for $|A|< 2\sqrt{p}$, (see for example \cite{Le87} for details) 
\[
N_p(H_p) \ll \sqrt{p}\log p (\log \log p)^2 .
\]
There is also a lower bound of the form $\gg \frac{\sqrt{p}}{\log p}$ for $|A| < \sqrt{p}$ with a bounded number of possible exceptions.

  Returning to methods from transcendence theory, Mit'kin \cite{Mi} and Heath-Brown \cite{He1} constructed two polynomials $F$ with $ d = p-1$ such that for all $a \mod p$ the corresponding $N_{p}(F - a)$ is bounded above by $Cp^{\frac{2}{3}}$ with estimates uniform in $a$. The polynomial considered by both authors was the truncated logarithm 
\begin{equation*}\label{eq0}
L(x) = x + \frac{x^2}{2} + \cdots + \frac{x^{p-1}}{p-1}.
\end{equation*}
 The truncated logarithm is not special in this regard and  Mit'kin established the same asymptotic bound for the truncated exponential $$E(x) = 1 + x + \frac{x^2}{2!} + \cdots + \frac{x^{p-1}}{(p-1)!}.$$ For the proof, one  constructs an auxiliary polynomial $\Phi$ (not vanishing identically) in several variables with relatively low degree in each variable, but with a high order of vanishing at each root of $F(x)$ in the prime field. Then, $N_{p}(F)$ is bounded by the degree of $\Phi$ divided by the order of vanishing. 

  To create $\Phi$ with high order of vanishing at the selected points, Leibniz' rule is used so that the higher derivatives of $\Phi$ are forced to  vanish at most of the roots of the original polynomial. In the case of the truncated logarithm and exponential, this is achieved using the fact that each satisfies a differential equation of the form 
\begin{equation}\label{eq1} 
\{x(1-x)\}^n f^{(n)}(x) \equiv a_n(x) + b_n(x)(x^p - x) + c_n(x)f(x) \mod p,
\end{equation}
for each $n=1,2,3,...$, where $a_n(x)$, $b_n(x)$ and $c_n(x)$ denote polynomials, of low degree relative to $n$, with integer coefficients (for the moment, we will be intentionally imprecise about what is meant by ``low degree"). Then since $f(x)$ is constant at our points of interest (and obviously as is $x^{p} -x$), the values taken by these (weighted) higher derivatives of $f(x)$ are those taken by certain polynomials of low degree.  It is then not difficult to determine that construction of $\Phi$ amounts to finding a nontrivial solution to a system of linear equations.  After this construction, one has to verify that when the variables in $\Phi$ are specialised for the problem under consideration, the resulting new polynomial, which is now of one variable, does not vanish identically (see Section 2 for some details).

It is interesting to note that $L$ is obtained by truncating a $G$-function, while $E$ is a truncation of an $E$-function (both types of functions were introduced by Siegel \cite{Siegel}; see the notes by Beukers \cite{Beukers} for recent results). Also the Hasse invariant $H_p$ above is a truncation of the hypergeometric function $_{2}F_{1}(\frac{1}{2},\frac{1}{2},1;x)$. It is then a natural question to ask if perhaps there is a much larger class of such polynomials $F$ for which there is a non-trivial estimate for $N_{p}(F)$. To this end, one could consider either $E$- or $G$-functions with rational Taylor coefficients but it is not clear how one should truncate these in a natural way. We illustrate this with the following example: let 
\[
r(x) = \sum_{k=0}^{\infty} \frac{2^{k}x^{2k +1}}{1\cdot3\cdot5\cdot ... \cdot(2k+1)},
\]
and let $R(x)$ be the polynomial obtained by truncating $r(x)$ with $0\leq k\leq \frac{p-3}{2}$. This polynomial satisfies a differential equation similar to \eqref{eq0} and we can show that $N_{p}(R) \ll p^{\frac{2}{3}}$ (we omit the details). Now consider instead the series 
\[
\sqrt{x}r(\sqrt{x})=\sum_{k=0}^{\infty} \frac{2^{k}x^{k +1}}{1\cdot3\cdot5\cdot ... \cdot(2k+1)}.
\] 
The polynomial associated with it should still be naturally truncated at $\frac{p-3}{2}$ (determined by the vanishing of the denominators of the coefficients)  but now the highest power of $x$ is too small so that we lack a formula of the type \eqref{eq0} involving $x^p -x$.

Our purpose in this note is extending the class of polynomials with high degree that have $o(p)$ roots as $p$ grows without bound, but which are obtained by truncating suitable functions that satisfy a higher order differential equation.  We observe that the methods of \cite{Mi} and \cite{He1} can be modified to accomplish this for truncations of polylogarithms and polyexponentials, for which we obtain an upper bound much weaker than a power saving in $p$. In general, the difficulty is twofold: the higher order derivatives depend on lower order derivatives (which are obviously non-constant at the roots of the polynomial) and the non-vanishing property requires, in essence, an algebraic independence involving the polynomial and its derivatives.

Our result, when applied to truncations of polylogarithms, may be stated as follows:

\begin{theorem1} Let $k \in \mathbb{N}$, and let $L_k$ denote the truncated polylogarithm $$L_k(x) = \sum_{i=1}^{p-1} \frac{x^i}{i^k}.$$ Then for $k\geq 2$ \begin{equation*} N_{p}(L_k) \ll_{k} \frac{p}{\log p}.\end{equation*} 
\end{theorem1}

Our analogous result for polyexponentials is similar:

\begin{theorem2} Let $k \in \mathbb{N}\cup \{0\}$, and let $E_k$ denote the truncated polyexponential $$E_k(x) = \sum_{i=1}^{p-1} \frac{x^i}{i! i^k}.$$ Then for $k\geq 1$ \begin{equation*} N_{p}(E_k) \ll_k \frac{p}{\sqrt{\log p}}. \end{equation*} \end{theorem2}

\begin{rem}The weakness in the results above are primarily due to the fact that we use an inductive method for $L_k$ and $E_k$. In employing a modification of the methods in \cite{Mi} and \cite{He1}, we are unable to control sufficiently the degrees of some polynomials that appear as coefficients involving $L_1,...,L_{k-1}$ and $E_0,...,E_{k-1}$, respectively so that some of the degrees grow exponentially.
\end{rem}

Finally, we consider the analogous question for the Bessel function
\[
J_0(x) = \sum_{n=0}^\infty \frac{(-1)^n\left(\frac{x}{2}\right)^{2n}}{(n!)^2},
\]
for which  we are able to save a power of $p$. Here we were unable to adapt the methods used above involving algebraic independence and so appeal to Siegel's original argument \cite{Siegel} showing the algebraic independence of $J_{0}(x)$ and $J'_{0}(x)$ over the complex numbers.  The argument can be applied with suitable modifications provided that the degree of the auxiliary polynomial in each variable is sufficiently small relative to $p$; this is a consequence of tame ramification, i.e., the ramification index at a point is relatively prime to $p$, which allows expansion of algebraic elements as a Puiseux series. This technique is not possible in general, as noted by Chevalley \cite{Ch}, as Puiseux series do not account for Artin-Schreier extensions.

For the truncation of the Bessel function, we assume that $p\geq 3$ and let 
\[
J_{0,p}(x) = \sum_{n=0}^{\frac{p + 1}{2}} \frac{(-1)^n \left(\frac{x}{2}\right)^{2n}}{(n!)^2}.
\]
Our result for the truncated Bessel function may then be stated as follows:

\begin{theorem3} For the Bessel function $J_0(x)$, one has $N_p(J_{0,p}) \ll p^{\frac{8}{9}}$. 
\end{theorem3}

The proof of Theorem 3 works equally well for the truncation of  Bessel functions of higher order
\[
J_\lambda(x) = \sum_{n=0}^\infty \frac{(-1)^n\left(\frac{x}{2}\right)^{2n+\lambda}}{n!(n+\lambda)!},
\]
provided that $\lambda$ is small relative to $p$.

\vskip 0.2in
\noindent{{\bf Acknowledgments.}} \\
Both authors thank Roger Heath-Brown for bringing their attention to the Hasse invariant and the corresponding comments.

\noindent AG thanks Alan Adolphson, Enrico Bombieri and Nick Katz for discussions over a period of time on the topics of this paper, and acknowledges support from the Vaughan Fund at OSU.

\noindent KW also thanks the Department of Mathematics at Oklahoma State University for their support during the writing of this paper. 
\vskip 0.2in
\section{Review of the Mit'kin/Heath-Brown construction.}

Here, we give the details of the basic mechanics of Stepanov's method as applied in \cite{Mi} and \cite{He1} to the case of $L(x)$ (the method for $E(x)$ is similar). One first constructs a polynomial $$\Phi(x,y,z) = \sum_{a,b,c} \lambda_{a,b,c} x^a y^b z^c \in \mathbb{F}_p[X,Y,Z]$$ with $$\deg_X \Phi \leq A, \; \deg_Y \Phi \leq B, \; \text{ and } \deg_Z \Phi \leq C;$$ and, furthermore, with the requirement that $\Psi(x)=\Phi(x,f(x),x^p)$ is not identically zero, but vanishes at each root of $L(x)$ in $\mathbb{F}_p$ with order at least $D$. This would guarantee that $$D N_0 \leq \deg \Psi(x) \leq A + (p-1)B + pC,$$ and thus yield a bound on $N_0$. By differentiating $\Psi(x)$ with use of \eqref{eq1} and observing that all terms of the form $x^p - x$ may be discarded, one finds that it is enough to require 
\begin{equation} \label{eq3} 
D(A+2D+C)< ABC
\end{equation} 
to guarantee the vanishing of $\Psi(x)$ at the zeros of $L(x)$ while maintaining that $\Phi(x,y,z)$ not vanish identically. The left-hand side of \eqref{eq3} simply reflects that we are arranging $\Psi^{(n)}(x)$ to vanish for $n < D$ at the roots of $L(x)$: the multiplier $D$ occurs because $\{\Psi^{(n)}(x)\}_{n=0}^{D-1}$ yields $D$ polynomials with coefficients that are linear forms in the coefficients of $\Phi$, and the term $A+2D+C$ represents a bound on the degree of the polynomials appearing in each higher derivative of $\Psi(x)$, which was obtained by use of \eqref{eq1}. On the right-hand side, the term $ABC$ is a consequence of the number of coefficients appearing in $\Phi$. 

The second part of the argument is deeper and reflects the transcendental nature of the power series giving rise to the polynomial.  One has to show that $\Psi(x)$ itself is not identically zero, and the arguments in \cite{Mi} and \cite{He1} are essentially the same on this, using the observation that $\Psi(x)$ will not vanish identically if it does not also vanish modulo $x^p$. This has the effect of putting restrictions on the parameters $A$, $B$ and $C$ above, namely that they cannot be chosen too small. Writing $$\Phi(x,y,z) = \sum_c F_c(x,y) z^c,$$ it is enough to show that the smallest $c$ in this sum for which $F_c(x,y)$ is nonzero satisfies $x^p \nmid F_c(x,L(x))$. This is Lemma 3 in \cite{He1}, and is established using Leibniz' rule, the binomial theorem, and an inductive argument. For $L(x)$, it is then enough to set $A = \lfloor p^{2/3} \rfloor$, $B = C = \lfloor p^{1/3} \rfloor$, and $D = \lfloor \frac{1}{3} p^{2/3} \rfloor$. 

\section{Proof of Theorem 1.\ }
\subsection{Construction of the auxiliary polynomial \ }\

We first establish a series of lemmas on $L_k(x)$. Our first lemma establishes that $L_k(x)$ satisfies a simple differential equation in terms of $L_1(x),...,L_{k-1}(x)$.

\begin{lemma1} The truncated $k$th polylogarithm $L_k(x)$ satisfies \begin{equation}\label{eq4} \{x(1-x)\}^n L_k^{(n)}(x) = a_{k,n}(x) + b_{k,n}(x)(x^p - x) + \sum_{i=1}^{k-1} c_{k,n,i}(x)L_i(x)\end{equation} for each $n=1,2,3,...$, where each $a_{k,n}(x)$, $b_{k,n}(x)$ and $c_{k,n,i}(x)$, with $i=1,..., k-1$ are polynomials with integer coefficients of degree at most $n+1$, $n-1$, and $n$, respectively. \end{lemma1}

\begin{proof}That this is true for $L_1(x) = L(x)$ is Lemma 1 of \cite{He1}. For $k>1$ and $n=1$, one has \begin{equation*} x(1-x)L_k'(x) = (1-x) L_{k-1}(x),\end{equation*} so one may set $a_{k,1}(x) = 0$, $b_{k,1}(x) = 0$, $c_{k,1,i}(x) = 0$ for $i=1,...,k-2$, and $c_{k,1,k-1}(x)=1-x$.

Assume that \eqref{eq4} has been verified for $k > 1$ and some $n \geq 1$. Differentiating \eqref{eq4} and multiplying by $x(x-1)$ yields a left-hand side equal to \begin{equation*} \{x(1-x)\}^{n+1} L_k^{(n+1)}(x)+ n(1-2x)\{x(1-x)\}^n L_k^{(n)}(x), \end{equation*} and for the right-hand side one obtains \begin{align*} & x(1-x) a_{k,n}'(x) + x(1-x)b_{k,n}'(x)(x^p - x) + \sum_{i=1}^{k-1} x(1-x)c_{k,n,i}'(x)L_i(x) \\& - x(1-x)b_{k,n}(x) + \sum_{i=1}^{k-2} (1-x)c_{k,n,i+1}(x)L_i(x) + c_{k,n,1}(x)(x^p - x). \end{align*} Therefore one may set \begin{equation*} a_{k,n+1}(x) = x(1-x)\{a_{k,n}'(x) - b_{k,n}(x)\} - n(1-2x)a_{k,n}(x),\end{equation*} \begin{equation*} b_{k,n+1}(x) = x(1-x)b_{k,n}'(x) + c_{k,n,1}(x)- n(1-2x)b_{k,n}(x),\end{equation*} \begin{equation*} c_{k,n+1,i}(x) = x(1-x)c_{k,n,i}'(x)+ (1-x)c_{k,n,i+1}(x) - n(1-2x)c_{k,n,i}(x) \end{equation*} for $i=1,...,k-2$, and \begin{equation*} c_{k,n+1,k-1}(x) = x(1-x)c_{k,n,k-1}'(x)-n(1-2x)c_{k,n,k-1}(x). \end{equation*} By the inductive argument, these polynomials possess integer coefficients and satisfy the required bounds on degrees. \end{proof}

Our second lemma, regarding a product of the functions $L_1(x),...,L_k(x)$, is essentially an application of Leibniz' rule, which allows us to bound the degree of the coefficients for terms appearing in higher derivatives. This will motivate our definition of the auxiliary polynomial.

\begin{lemma2} Let $a,c_1,...,c_k$ be nonnegative integers. The derivative $$\{x(1-x)\}^n\frac{d^n}{dx^n}\{x^{a}L_1(x)^{c_1} \cdots L_k(x)^{c_k}\},$$ if not identically zero, is equal to a sum of polynomials of the form \begin{equation*} p(x) q_1(x) \cdots q_l(x) L_1(x)^{d_1} \cdots L_k(x)^{d_k}\end{equation*} modulo $x^p - x$ with integer coefficients, where \begin{equation*} \sum_{i=1}^k id_i \leq \sum_{i=1}^k ic_i, \end{equation*} \begin{equation*} 0 \leq l \leq \min\left(\sum_{i=1}^k c_k,n\right), \end{equation*} and $p(x),q_1(x),...,q_l(x)$ are polynomials where \begin{equation*}\deg q_{i}(x) \leq g_i + 1,\; i=1,...,l,\end{equation*} and \begin{equation*} \deg p(x) = a + n - \sum_{i=1}^l g_i. \end{equation*} \end{lemma2}

In the proof of Theorem 1, we will need to use the fact that the auxiliary polynomial we define does not vanish identically, and this is precisely why the sum $\sum_{i=1}^k ic_i$ appears. In the case of $k=1$, this is unimportant, as powers of $L_1(x)$ simply decrease. For general $k$, the same is not necessarily true: as in Lemma 1, powers of $L_j(x)$ ($j \leq k$) are exchanged for powers of $L_1(x),...,L_{j-1}(x)$.

\begin{proof} The case of $k=1$ follows as in \cite{He1}: By Lemma 1, one may write \begin{equation*} \{x(1-x)\}^m L_1^{(m)}(x)=a_{1,m}(x) + b_{1,m}(x)(x^p - x)\end{equation*} where the degrees of $a_{1,m}(x)$ and $b_{1,m}(x)$ are less than $m+1$ and $m-1$, respectively, and all coefficients are integers. Successive application of this property implies that $$\{x(1-x)\}^{n-r}\frac{d^{n-r}}{dx^{n-r}}(x^{a}) \{x(1-x)\}^{r}\frac{d^{r}}{dx^{r}}(L_1(x)^{c_1})$$ either vanishes or is a sum of polynomials of the form \begin{equation*} p(x) q_1(x) \cdots q_l(x) L_1(x)^{c_1 - l}\end{equation*} modulo $x^p - x$ with integer coefficients, where $0 \leq l \leq \min \{c_1,n\}$ and $p(x),q_1(x),...,q_l(x)$ are certain polynomials. We let each of the polynomials $q_1(x),...,q_l(x)$ equal $a_{1,m}(x)$ for some $m$ by Leibniz' rule, and therefore the sum of their degrees is at most $r+l$; in fact, we may thus define $g_1,...,g_l$ so that $\deg q_i \leq g_i + 1$ for each $i=1,...,l$ and \begin{equation*} \sum_{i=1}^l (g_i + 1) = \left(\sum_{i=1}^l g_i \right) + l \leq r + l. \end{equation*} The polynomial $p(x)$ will simply represent $\{x(1-x)\}^{n-r}\frac{d^{n-r}}{dx^{n-r}}(x^{a})$, and is thus a polynomial of degree $n-r+a$. As \begin{equation*}\{x(1-x)\}^n\frac{d^n}{dx^n}\{x^{a}L_1(x)^{c_1}\}=\sum_{r=0}^n \binom{n}{r} \{x(1-x)\}^{n-r}\frac{d^{n-r}}{dx^{n-r}}(x^{a})  \{x(1-x)\}^{r}\frac{d^{r}}{dx^{r}}(L_1(x)^{c_1}), \end{equation*} the result then follows for the case of $k=1$. 

If $k > 1$, a similar argument applies: Setting $r_1,...,r_k$ nonnegative with $r = r_1 + \cdots r_k$, Lemma 1 yields that \begin{equation*} \{x(1-x)\}^{n-r}\frac{d^{n-r}}{dx^{n-r}}(x^{a}) \{x(1-x)\}^{r_1}\frac{d^{r_1}}{dx^{r_1}}(L_1(x)^{c_1})\cdots \{x(1-x)\}^{r_k}\frac{d^{r_k}}{dx^{r_k}}(L_1(x)^{c_1}) \end{equation*} may be written as a sum of polynomials of the form \begin{equation*} p(x) q_1(x),...,q_l(x) L_1(x)^{d_1} \cdots L_k(x)^{d_k} \end{equation*} modulo $x^p - x$ with integer coefficients, where $p(x),q_1(x),...,q_l(x)$ are certain polynomials. In this case, \eqref{eq4} has been applied to derivatives of $L_1(x),...,L_k(x)$ according to Leibniz' rule and gives $\sum_{i=1}^k id_i \leq \sum_{i=1}^k ic_i$. As in the case of $k=1$, $p(x)$ represents $\{x(1-x)\}^{n-r}\frac{d^{n-r}}{dx^{n-r}}(x^{a})$, and either vanishes or is a polynomial of degree $n-r+a$. Furthermore, Lemma 1 implies that the polynomials $q_i(x)$ may have degree bounded by $g_i$ rather than $g_i+1$, depending on whether $q_i(x)$ occurs as a coefficient of a polylogarithm; the bounds occurring in the case of $k=1$ are thus again valid, and the result follows. \end{proof}

Finally, in order to prove the nonvanishing of the auxiliary polynomial, we will need to prove that at least one of its coefficients does not vanish. Our third lemma identifies a nonvanishing term in a higher derivative of a product of the functions $L_1(x),...,L_k(x)$. We introduce the following notation: For a polynomial $f$, possibly of multiple variables, let $a(f)$ denote the degree of $f$ in its first variable. If $f$ is a polynomial in $k+1$ variables $z,x_1,...,x_k$, let $\tilde{f}(x):=f(x,L_1(x),...,L_k(x))$ and $Sf(x):= \{x(1-x)\}^{a(f) + 1} \tilde{f}^{(a(f)+1)}(x)$; furthermore, let $c(f)$ denote the largest sum $\sum_{i=1}^k ic_i$ of powers $(c_1,...,c_k)$ appearing as a product $g(z) x_1^{c_1} \cdots x_k^{c_k}$ in the expression for $f$.

Let us now define a class of functions that will be useful for the proof of Theorem 1. 

\begin{definition} Let $\mathcal{S}$ denote the class of polynomials $f$ in $k+1$ variables with coefficients in $\mathbb{F}_p$ and $\max\{a(f),1\} \cdot 4^{c(f)} < p -1$. \end{definition}

Such a polynomial $f$ will appear in our auxiliary polynomial, which will not vanish $\mod x^p$.

\begin{lemma3} Suppose that \emph{(1)} $f \in \mathcal{S}$ and \emph{(2)} $0 < \max_i \deg_{x_i} f < p$. Then $$Sf(x) = G_f(x,L_1(x),...,L_k(x))\mod x^p$$ where $G_f \in \mathcal{S}$ does not vanish identically; in particular, $x^p \nmid \tilde{f}(x)$.\end{lemma3}

\begin{proof} Let $x_1^{c_1} \cdots x_k^{c_k}$ be the product in $f$ for which (1) the quantity \begin{equation*}2c_1 + 3 c_2 + \cdots + (k+1)c_k \end{equation*} is largest, and (2) for all other terms $x_1^{c_1'} \cdots x_k^{c_k'}$ in $f$ with $\sum_{i=1}^k (i+1)c_i = \sum_{i=1}^k (i+1)c_i'$, it holds that $c_r' > c_r$ implies $c_s > c_s'$ for some $s > r$. Let $g(z)$ denote the coefficient of $x_1^{c_1} \cdots x_k^{c_k}$. 

If $c_1 \neq 0$, then by choice of $g(z)x_1^{c_1} \cdots x_k^{c_k}$ it follows from Lemma 1 that the term in $Sf(x)$ containing the exact product $L_1(x)^{c_1 - 1} L_2(x)^{c_2} \cdots L_k(x)^{c_k}$ is obtained from only \begin{equation}\label{eq7} (1-x)^{a(f) + 1} \frac{d^{a(f) + 1}}{dx^{a(f) + 1}} (g(x)L_1(x)^{c_1})L_2(x)^{c_2} \cdots L_k(x)^{c_k}.\end{equation} By the proof of Lemma 5 of \cite{He1}, the component of the coefficient of \begin{equation}\label{eq8} x^{a(g)}L_1(x)^{c_1 - 1} L_2(x)^{c_2} \cdots L_k(x)^{c_k}\end{equation} in $Sf(x)$ obtained  from \eqref{eq7} is nonzero. As this is the only contribution to the coefficient of \eqref{eq8}, it follows that the coefficient of $$x^{a(g)+a(f)+1}L_1(x)^{c_1 - 1} L_2(x)^{c_2} \cdots L_k(x)^{c_k}$$ in $Sf(x)$ is nonzero. 

Suppose then that $c_1 = c_2 = ... = c_{j-1} = 0$ and $c_j \neq 0$. The function $L_j(x)$ satisfies \begin{equation}\label{eq9} x^{l}L_j^{(l)}(x) = (-1)^{l-1}(l-1)!L_{j-1}(x) + g_{j,l}(x), \end{equation} where, as in Lemma 1, the function $g_{j,l}(x)$ is a linear combination of $L_1(x)$,...,$L_{j-2}(x)$, $x^p-x$, and $1$, with coefficients equal to polynomials of low degree in $x$. As in the previous case, the term in $Sf(x)$ containing the exact product \begin{equation*} L_{j-1}(x)L_j(x)^{c_j - 1}L_{j+1}(x)^{c_{j+1}}\cdots L_k(x)^{c_k} \end{equation*} is obtained from only \begin{equation*} x^{a(f) + 1} \frac{d^{a(f) + 1}}{dx^{a(f) + 1}}(g(x)L_j(x)^{c_j})L_{j+1}(x)^{c_{j+1}}\cdots L_k(x)^{c_k} .\end{equation*} By Leibniz' rule, we may write \begin{align}\label{eq11} &\frac{d^{a(f) + 1}}{dx^{a(f) + 1}}(x^{a(g)}L_j(x)^{c_j}) \\&\notag= \sum_{l=0}^{a(f) + 1} \binom{a(f) + 1}{l} \frac{d^l}{dx^l}(x^{a(g)})\frac{d^{a(f) + 1 - l}}{dx^{a(f) + 1 - l}}(L_j(x)^{c_j}). \end{align} Furthermore, we have \begin{align}\label{eq12}& \frac{d^{a(f) + 1 - l}}{dx^{a(f) + 1 - l}}(L_j(x)^{c_j}) \\&\notag= \sum_{l_1,...,l_{c_j}} \binom{a(f) + 1 - l}{l_1,...,l_{c_j}} L_j^{(l_1)}(x) \cdots L_j^{(l_{c_j})}(x), \end{align} where $\sum_{r=1}^{c_j} l_r = a(f) + 1 - l$. It follows from \eqref{eq9}, \eqref{eq11} and \eqref{eq12} that the coefficient of \\$x^{a(g)}L_{j-1}(x)L_j(x)^{c_j - 1}L_{j+1}(x)^{c_{j+1}}\cdots L_k(x)^{c_k}$ in $Sf(x)$ is equal to $$c_j \sum_{l=0}^a \binom{a(f) + 1}{l} \frac{a(g)!}{(a(g) - l)!} (-1)^{a(f) - l}(a(f) - l)!.$$ As $c_j < p$, it follows as in the proof of Lemma 5 of \cite{He1} that this sum is nonzero in $\mathbb{F}_p$. Again, the coefficient of $x^{a(g)+a(f)+1}L_1(x)^{c_1 - 1} L_2(x)^{c_2} \cdots L_k(x)^{c_k}$ in $Sf(x)$ is nonzero. 

The existence of $G_f$ has thus been established. By Lemma 1 and Leibniz' rule, it follows that that $a\left(G_f\right) \leq 3a(f) +2$ and $c\left(G_f\right) \leq c(f) - 1$. Therefore \begin{align*}  \max  \left\{a\left(G_f\right),1\right\}\cdot 4^{c\left(G_f\right)} &\leq \max\{3a(f) +2,1\}\cdot 4^{c(f) - 1}\\& = \max\left\{\frac{1}{4}(3a(f) +2),\frac{1}{4}\right\}\cdot 4^{c(f)}\\& \leq \max\{a(f),1\} \cdot 4^{c(f)} < p - 1. \end{align*} Suppose that $c(f) = 0$; then it is obvious by definition of $\mathcal{S}$ that $x^p \nmid \tilde{f}(x)$. If $c(f) > 0$, suppose that $x^p | \tilde{f}(x)$; then $x^p | \tilde{f}^{(a(f)+1)}(x)$, and thus $x^p | Sf(x)$. It follows that $$G_f(x,L_1(x),...,L_k(x)) \equiv Sf(x) \equiv 0 \mod x^p,$$ which contradicts the induction hypothesis. \end{proof}

We are now prepared to present the proof of Theorem 1.

\begin{proof}[Proof of Theorem 1] Let us define 
\begin{equation*} \Phi_k(x_{-1},x_0,x_1,...,x_k) = \sum_{a,b,c_1,...,c_k \geq 0} \lambda_{a,b,c_1,...,c_k} x_{-1}^{a} x_0^{b} x_1^{c_1}\cdots x_k^{c_k}. \end{equation*}
Put
\begin{equation*}\Psi_k(x) = \Phi_k(x,x^p,L_1(x),L_2(x),...,L_k(x)),\end{equation*} 
so that \begin{equation*} \Psi_k(x) = \sum_{a,b,c_1,...,c_k \geq 0} \lambda_{a,b,c_1,...,c_k} x^{a} x^{pb} L_1(x)^{c_1} \cdots L_k(x)^{c_k}. \end{equation*}
Our goal is to retrieve an expression for the higher derivatives of $\Psi_k(x)$ modulo $x^p - x$ in terms of $L_1(x)$, $L_2(x)$,..., $L_{k-1}(x)$, with coefficients equal to polynomials of low degree in $x$.

Let \begin{equation*} \deg_{x_{-1}} \Phi_k < A \text{ and } \deg_{x_0} \Phi_k < B. \end{equation*} Furthermore, let us require that the largest sum $\sum_{i=1}^k ic_i$ of $(c_1,...,c_k)$ with $c_i \geq 0$ appearing together as a product of powers of $x_1,...,x_k$ in $\Phi_k$ is at most $C$. We wish to select the coefficients of $\Phi_k$ appropriately to guarantee that $\Psi_k^{(n)}(x)$ vanishes at almost all zeros of $L_k(x)$ for all $n < D$, with  $D$ to be chosen (\emph{caveat lector}: the labelling used here is somewhat different from that appearing in \cite{He1}). By appropriate selection of $A$, $B$, $C$, and $D$, it will suffice to require that \begin{enumerate}[(i)] \item $\{x(1-x)\}^{n}\Psi_k^{(n)}(x)|_{x=\alpha} = 0$ for each $n < D$ and zero $\alpha \in \mathbb{F}_p$ of $L_k(x)$; and 
\item  $\Psi_k$ does not vanish identically as a polynomial. \end{enumerate}

Let $S(C,k)$ denote the number of possible values of $(c_1,...,c_k)$. The function $\Phi_k$ will thus have $AB \cdot S(C,k)$ unknowns $\lambda_{a,b,c_1,...,c_k}$ that we must select, as in the right-hand side of \eqref{eq3}. A term of the form $x^{a} x^{pb} L_1(x)^{c_1} \cdots L_k(x)^{c_k}$ appearing in the expression of $\Psi_k$ satisfies \begin{equation*}  \frac{d^n}{dx^n}\{x^{a} x^{pb} L_1(x)^{c_1} \cdots L_k(x)^{c_k}\} = x^{bp} \frac{d^n}{dx^n}\{x^{a}L_1(x)^{c_1} \cdots L_k(x)^{c_k}\}. \end{equation*} The polynomial $\{x(1-x)\}^n \frac{d^n}{dx^n}(x^{a})$ is either identically zero or of degree equal to $a+n$, and we have $x^{bp} \equiv x^b \; \text{mod}\; x^p - x$. By Lemma 2, we may therefore write \begin{align}\label{eq10}\notag \{x(1&-x)\}^{n} \Psi_k^{(n)}(x) \\&\equiv \sum_{\substack{ d_1,...,d_k \\ d_1 + \cdots +  d_k < C}} f(x;k,d_1,...,d_k,n)L_1(x)^{d_1} \cdots L_k(x)^{d_k}\;\;\mod\; x^p - x \end{align} in $\mathbb{F}_p[x]$, where $\deg_x f(x;k,d_1,...,d_k,n) < A + B + 2n$ for each $d_1,...,d_k$ and $n < D$. As we are considering only the zeros of $L_k(x)$ in $\mathbb{F}_p$, we may disregard all terms in \eqref{eq10} where $d_k$ is nonzero, as well as any terms in $\{x(1-x)\}^{n}\Psi_k^{(n)}$ where $x^p - x$ appears. Therefore our system of coefficients for $\Psi_k^{(n)}$ is in number at most $(A+B + 2D)\cdot S(C,k-1)$. Also, the coefficients of $\Psi_k^{(n)}$ are linear forms in the coefficients of $\Phi_k$. As we are considering all $n < D$, there are $D$ such systems of coefficients. Provided that  \begin{equation*} D(A+B + 2D) \cdot S(C,k-1) < AB \cdot S(C,k), \end{equation*} there will exist a nontrivial choice of coefficients of $\Phi_k$ for which $\Psi_k^{(n)}$ vanishes at the zeros of $L_k(x)$ for $n < D$ (excepting $0$ and $1$). This concludes the proof of (i).

For the proof of (ii), we must verify that $\Psi_k(x)$ does not vanish identically with this choice of coefficients. We may write \begin{equation*} \Phi_k(x_{-1},x_0,x_1,...,x_k) = \sum_b f_b(x_{-1},x_1,...,x_k) x_0^b. \end{equation*} Let $b_0$ be the smallest value of $b$ so that $f_b(x_{-1},x_1,...,x_k) \neq 0$. Such a $b_0$ exists because $\Phi_k$ is not identically zero by selection of coefficients in the first step of the proof. If $\Psi_k$ were identically zero, then $\tilde{f}_{b_0}(x)$ would be divisible by $x^p$, which is not possible by Lemma 3 if $f_{b_0} \in \mathcal{S}$. It is enough, then, for us to choose $A$, $B$, $C$, and $D$ to satisfy \begin{enumerate} \item $D(A+B + 2D) \cdot S(C,k-1) < AB \cdot S(C,k)$, and \item $A \cdot 4^C < p -1$. \end{enumerate} 
We let $A = D =(\log p)^2$ and $C = \varepsilon \log p$, where $\varepsilon$ is chosen suitably small so that condition (2) is satisfied. As $S(C,k) = \frac{C^k}{(k!)^2} + O(C^{k-1}),$ we may let $B = R \log p$ for sufficiently large $R$. We then have $$N_{p}(L_{k}) \ll(A + pB + (p-1)C)/D \ll p/\log p .$$  This latter follows from the discussion in Section 2 and using the fact that the contributions of $L_1(x),...,L_k(x)$ to the degree of $\Psi_k$ as a polynomial in a single variable appear as products $L_1(x)^{c_1} \cdots L_k(x)^{c_k}$ that satisfy $$\sum_{i=1}^k c_i \leq \sum_{i=1}^k ic_i \leq C.$$    \end{proof}

\section{Proof of Theorem 2.\ }

As with $L_k(x)$, we require a few preliminary results on $E_k(x)$.

\begin{lemma4} The truncated $k$th polyexponential $E_k(x)$ satisfies \begin{equation}\label{eq13} x^n E_k^{(n)}(x) = a_{k,n}(x) + b_{k,n}(x)(x^p - x) + c_{k,n,0}(x)E_0(x)+ \sum_{i=1}^{k-1} c_{k,n,i}(x)E_i(x)\end{equation} for each $n=1,2,3,...$, where each $a_{k,n}(x)$, $b_{k,n}(x)$ and $c_{k,n,i}(x)$, with $i=0,..., k-1$ are polynomials with integer coefficients of degree at most $n$, $n-1$, and $n$, respectively. \end{lemma4}

\begin{proof} The case of $E_0(x) = E(x)$ is Lemma 2 of \cite{Mi}. If $k >1$, we have \begin{equation*} xE_k'(x) = E_{k-1}(x), \end{equation*} so one may set $a_{k,1}(x) = 0$, $b_{k,1}(x) = 0$, $c_{k,1,i}(x) = 0$ for $i=0,...,k-2$, and $c_{k,1,k-1}(x) = 1$. For the inductive step, differentiating \eqref{eq13} and multiplying by $x$ gives a left-hand side equal to \begin{equation*} nx^n E_k^{(n)}(x) + x^{n+1}E_k^{(n+1)}(x). \end{equation*} For the right-hand side, we obtain \begin{align*} & xa_{k,n}'(x) + xb_{k,n}'(x)(x^p - x) - xb_{k,n}(x) + xc_{k,n,0}'(x)E_0(x) \\& + c_{k,n,0}(x)(xE_0(x) + (x^p - x) + x) + \sum_{i=1}^{k-1} xc_{k,n,i}'(x)E_i(x) + \sum_{i=0}^{k-2} c_{k,n,i+1}(x)(E_i(x) - 1). \end{align*} Therefore one may set \begin{equation*} a_{k,n+1}(x) = xa_{k,n}'(x)+xb_{k,n}(x)+xc_{k,n,0}(x) + \sum_{i=0}^{k-2} c_{k,n,i+1}(x) - na_{k,n}(x),\end{equation*} \begin{equation*} b_{k,n+1}(x) = xb_{k,n}'(x) + c_{k,n,0}(x)-nb_{k,n}(x),\end{equation*} \begin{equation*} c_{k,n+1,0}(x) = xc_{k,n,0}(x) + c_{k,n,1}(x) - nc_{k,n,0}(x), \end{equation*} \begin{equation*} c_{k,n+1,i}(x) = xc_{k,n,i}'(x) + c_{k,n,i+1}(x)-nc_{k,n,i}(x) \end{equation*}for $i=1,...,k-2$, and \begin{equation*} c_{k,n+1,k-1}(x) = xc_{k,n,k-1}'(x)- nc_{k,n,i}(x). \end{equation*} By the inductive argument, these polynomials possess integer coefficients and satisfy the required bounds on degrees. \end{proof}

\begin{lemma5} Let $a,c_0,...,c_k$ be nonnegative integers. The derivative $$x^n\frac{d^n}{dx^n}\{x^{a}E_0(x)^{c_0} \cdots E_k(x)^{c_k}\},$$ if not identically zero, is equal to a sum of polynomials of the form \begin{equation*} p(x) q_1(x) \cdots q_l(x) E_0(x)^{d_0} \cdots E_k(x)^{d_k}\end{equation*} modulo $x^p - x$ with integer coefficients, where \begin{equation*} \sum_{i=1}^k id_i \leq \sum_{i=1}^k ic_i, \end{equation*} \begin{equation*} 0 \leq l \leq \min\left(\sum_{i=0}^k c_k,n\right), \end{equation*} and $p(x),q_1(x),...,q_l(x)$ are polynomials where \begin{equation*}\deg q_{i}(x) \leq g_i,\; i=1,...,l,\end{equation*} and \begin{equation*} \deg p(x) = a + n - \sum_{i=1}^l g_i. \end{equation*} \end{lemma5}

\begin{proof} This follows from Lemma 4; the proof is similar to that of Lemma 2. \end{proof}

If $f$ is a polynomial in $k+2$ variables $z,x_0,...,x_k$, let $$\tilde{f}(x) := f(x, E_0(x), E_1(x),..., E_k(x))$$ and $Tf(x) := x^{a(f) + 1} \tilde{f}^{(a(f)+1)}(x)$. Also, let $d_i(f)$ denote the degree of $f$ in $x_i$ for each $i=0,1,...,k$, and let $c(f)$ denote the largest sum $\sum_{i=1}^k ic_i$ of powers $(c_1,...,c_k)$ appearing as a product $g(z,x_0) x_1^{c_1} \cdots x_k^{c_k}$ in the expression for $f$.

Let us now define a class of functions for $E_k(x)$ analogous to the class $\mathcal{S}$ for $L_k(x)$.  

\begin{definition} Let $\mathcal{T}$ denote the class of polynomials $f$ in $k+1$ variables with coefficients in $\mathbb{F}_p$ and $(a(f)+1)\cdot (d_0(f)+a(f)+2) \cdot 6^{c(f)} < p -1$. \end{definition}

\begin{lemma6} Suppose that \emph{(1)} $f \in \mathcal{T}$ and \emph{(2)} $0 < \max_i \deg_{x_i} f < p$. Then $$Tf(x) = H_f(x,E_0(x),E_1(x),...,E_k(x))\text{ mod }x^p$$ where $H_f \in \mathcal{T}$ does not vanish identically; in particular, $x^p \nmid \tilde{f}(x)$. \end{lemma6}

\begin{proof} Suppose first that none of $x_1,...,x_k$ appear in the expression for $f$. In this case we may write $f = \sum_{i=0}^{d_1} g_i(z)x_0^i$. If $Tf(x) \equiv 0 \text{ mod } x^p$, then $x^{p - (a(f) + 1)} \mid \tilde{f}^{(a(f)+1)}$(x). By the proof of Lemma 4 of \cite{Mi}, this is impossible. 

Otherwise, let $g(z,x_0)x_1^{c_1} \cdots x_k^{c_k}$ be the component of $f$ for which (1) the quantity $\sum_{i=1}^k (i+1) c_i$ is largest, and (2) for all other terms $h(z,x_0)x_1^{c_1'} \cdots x_k^{c_k'}$ in $f$ with $\sum_{i=1}^k (i+1)c_i = \sum_{i=1}^k (i+1)c_i'$, it holds that $c_r' > c_r$ $(r \geq 1)$ implies $c_s > c_s'$ for some $s > r$ $(s \geq 1)$. If $c_1 \neq 0$, then by choice of $g(z,x_0)x_1^{c_1} \cdots x_k^{c_k}$ it follows from Lemma 4 and $x^n E_0^{(n)}(x) \equiv x^n E_0(x)\text{ mod } x^p$ $(n \geq 1)$ that the term in $Tf(x)$ containing the exact product $E_1(x)^{c_1 - 1} E_2(x)^{c_2} \cdots E_k(x)^{c_k}$ is obtained from only \begin{equation}\label{eq14}x^{a(f) + 1} \frac{d^{a(f) + 1}}{dx^{a(f) + 1}} (g(x,E_0(x))E_1(x)^{c_1})E_2(x)^{c_2} \cdots E_k(x)^{c_k}.\end{equation} Let $g(x,E_0(x)) = \sum_{i=0}^{c_0} g_i(x)E_0(x)^i$, and let $a$ be the degree of $g_{c_0}(x)$ in $x$. As $xE_1'(x)= E_0(x)-1$ and $x^{a(f)} E_0^{a(f)}(x) \equiv x^{a(f)} E_0(x)\text{ mod } x^p$, the coefficient of $$x^{a+a(f)}E_0(x)^{c_0+1}E_1(x)^{c_1 - 1} E_2(x)^{c_2} \cdots E_k(x)^{c_k}$$ obtained from \eqref{eq14} is equal to $$\sum_{k_1,...,k_{c_0}} \binom{a(f)}{k_1,...,k_{c_0}} = c_0^{a(f)} \not\equiv 0 \text{ mod } p.$$

Suppose then that $c_1 = c_2 = ... = c_{j-1} = 0$ and $c_j \neq 0$. As in the previous case, the term in $Tf(x)$ containing the exact product \begin{equation*} \label{eq15} E_{j-1}(x)E_j(x)^{c_j - 1}E_{j+1}(x)^{c_{j+1}}\cdots E_k(x)^{c_k} \end{equation*} is obtained from only \begin{equation} x^{a(f) + 1} \frac{d^{a(f) + 1}}{dx^{a(f) + 1}}(g(x,E_0(x))E_j(x)^{c_j})E_{j+1}(x)^{c_{j+1}}\cdots E_k(x)^{c_k}.\end{equation} Once again, as $xE_j'(x)= E_{j-1}(x)-1$ and $x^{a(f)} E_0^{a(f)}(x) \equiv x^{a(f)} E_0(x)\text{ mod } x^p$, the coefficient of $$x^{a+a(f)}E_0(x)^{c_0}E_{j-1}(x)E_j(x)^{c_j-1}E_{j+1}(x)^{c_{j+1}}\cdots E_k(x)^{c_k}$$ obtained from \eqref{eq15} is nonzero.  

The existence of $H_f$ has thus been established. By Lemma 1 and Leibniz' rule, it follows that that $a(H_f) \leq 2a(f) +1$, $d_0(H_f) \leq d_0(f) + a(f) + 1$, and $c\left(H_f\right) \leq c(f) - 1$. Therefore \begin{align*}  (a(H_f) + 1)\cdot(d_0(H_f)+a(H_f)+2) \cdot 6^{c\left(H_f\right)} &\leq (2a(f) + 2) \cdot (d_0(f) + 3a(f)+4) \cdot 6^{c(f) - 1} \\& =(a(f) + 1) \cdot \frac{1}{3}(d_0(f)+ 3a(f)+4) \cdot 6^{c(f)} \\& \leq (a(f) + 1)\cdot (d_0(f)+a(f)+2) \cdot 6^{c(f)} \\&< p-1.\end{align*} The remainder of the proof follows as in Lemma 3. \end{proof}

\begin{proof}[Proof of Theorem 2] As with $L_k(x)$, it is necessary to construct an auxiliary polynomial $\Phi_k$, but the proof mirrors that of Theorem 1. In fact, with $C$ as the bound on $\sum_{i=1}^k ic_i$ and $E$ as the bound on the degree in $E_0(x)$ for $\Phi_k$, our desired bounds are \begin{enumerate} \item $D(A+B +D) \cdot (E+D) \cdot S(C,k-1) < ABE \cdot S(C,k)$, and \item $(A+1) \cdot (E+A+2) \cdot 6^C < p -1$, \end{enumerate} where (2) is necessary to account for the fact that $E_0(x)$ does not vanish in its derivatives. We let $A=D=(\log p)^{2}$, $E = (\log p)^{\frac{3}{2}}$, and $C = \varepsilon \log p$, where $\varepsilon$ is chosen suitably small so that condition (2) is satisfied. As $S(C,k) = \frac{C^k}{(k!)^2} + O(C^{k-1}),$ we may let $B = R (\log p)^{\frac{3}{2}}$ for sufficiently large $R$. We then have $$N_{p}(E_{k}) \ll(A + pB + (p-1)(E+C)/D \ll p/\sqrt{\log p} .$$\end{proof}

\section{Proof of Theorem 3.\ }

The truncated Bessel function $J_{0,p}(x)$ satisfies the differential equation 
\begin{equation}\label{eqA1}
J_{0,p}''(x)+\frac{1}{x}J_{0,p}'(x) + J_{0,p}(x) \equiv 0 \mod x^p,
\end{equation}
and furthermore satisfies
\[
xJ_{0,p}''(x)+J_{0,p}'(x) + xJ_{0,p}(x) = \frac{(-1)^{\frac{p+1}{2}}}{2^{p+1}\left(\left(\frac{p+1}{2}\right)!\right)^2}x^2(x^p - x) + \frac{(-1)^{\frac{p+1}{2}}}{2^{p+1}\left(\left(\frac{p+1}{2}\right)!\right)^2}x^3.
\]

For the proof of Theorem 3, we require a preliminary lemma, which establishes a form of transcendence for the truncated Bessel function.

\begin{lemma7} Suppose that $n \in \mathbb{N}$ with $n^2 < p$. The function $J_{0,p}(x)$ is not a solution to any nonzero equation \begin{equation} \label{eq16} a_n(x)T^n + a_{n-1}(x)T^{n-1} + \cdots + a_0(x) \equiv 0 \mod x^p \end{equation} with $a_0(x),...,a_n(x) \in \overline{\mathbb{F}}_p[x]$ and $\max_i \deg a_i(x) < n$. \end{lemma7}

\begin{proof} Suppose that $y = J_{0,p}(x)$ is a solution to an equation \eqref{eq16} as in the statement of the Lemma. Let $z \in \overline{\mathbb{F}_p(x)}$ be a solution to the equation \begin{equation*} f(T)= a_n(x)T^n + a_{n-1}(x)T^{n-1} + \cdots + a_0(x)= 0. \end{equation*} As $n < p$, it follows that the extension $\overline{\mathbb{F}}_p(x)(z)|\overline{\mathbb{F}}_p(x)$ is tamely ramified. Thus (see, for example, \cite{Ch}) $z$ admits an expression at $x=0$ of the form $$z = \sum_{k=0}^\infty b_k x^{r_k}$$ with $r_0 < r_1 < \cdots$ rational exponents. Therefore the expression $y - z = g(x) \;\text{mod}\;x^p$ is a well-defined Puiseux series. Furthermore with $\alpha_1,...,\alpha_n$ the roots of $f(T)$, we have \begin{equation*} \prod_{i=1}^n (y - \alpha_i) \equiv 0 \mod x^p, \end{equation*} from which we conclude that for some $z = \alpha_i$, \begin{equation*} y - z \equiv 0 \mod x^{\lfloor \frac{p}{n}\rfloor}.\end{equation*} It follows from \eqref{eqA1} that \begin{equation} \label{eq999} x^2 z'' +xz' + x^2 z \equiv x^2 y'' +xy' + x^2 y \equiv 0 \mod x^{\lfloor \frac{p}{n} \rfloor}. \end{equation} As $\max_i \deg a_i(x) < n$ and the ramification index of any point is bounded by the degree of the extension, it follows that the degree of $z$ at any point, whether as a pole or zero, cannot be greater than $n^2$. Thus we may write an expression for $z$ at a branch of infinity as $$z = \sum_{k=0}^\infty c_k x^{s_k}$$ with $s_0 > s_1 > \cdots$ rational exponents, $-n^2 \leq s_0 \leq n^2$, and $c_0 \neq 0$. As $n^3 < p$, application of \eqref{eq999} yields that $c_0 = 0$, a contradiction. \end{proof}

We are now prepared to give a proof of Theorem 3. 

\begin{proof}[Proof of Theorem 3] The first step proceeds as in the proofs of Theorems 1 and 2, with the construction of an auxiliary polynomial, which in this case is a function of four variables. With $y = J_{0,p}(x)$, we set $\Psi(x) = \Phi(x,x^p,y,y'),$ where in this case we require that the total degree of $\Phi$ as a function of its third and fourth variables be at most $C$. As the number of nonnegative integer solutions to the inequality $x_1 + x_2 \leq C$ is simply $\frac{(C+1)(C+2)}{2}$, we obtain a now familiar bound: \begin{equation*} D(A + B + 2D) < AB\cdot \frac{(C+2)}{2}.\end{equation*}

For the second step of the proof, we suppose that there exists a nonzero polynomial $P(x_1,x_2)$ with coefficients in $\overline{\mathbb{F}}_p[x]$ of degree at most $s$ and total degree in $x_1$ and $x_2$ at most $t$ that satisfies \begin{equation} \label{eq20} P(y,y') \equiv 0 \mod x^p. \end{equation} Also, let \begin{equation} \label{eq21} P^*(x_1,x_2) =  P_x(x_1,x_2) + x_2 P_{x_1}(x_1,x_2) - (x_1 + \frac{1}{x}x_2)P_{x_2}(x_1,x_2).\end{equation} By the differential equation \eqref{eqA1} for $y$, we have by \eqref{eq21} that \begin{equation*} \frac{d}{dx}P(y,y') \equiv P^*(y,y') \mod x^p, \end{equation*} from which it follows that \begin{equation} \label{eq23} P^*(y,y') \equiv 0 \mod x^p.\end{equation} Let $R(y)$ denote the resultant of $P(y,x_2)$ and $xP^*(y,x_2)$ as polynomials in the second variable. By \eqref{eq20} and \eqref{eq23}, we have for suitable polynomials $u$ and $v$ (in $x_2$ and the coefficients of $P(y,x_2)$ and $xP^*(y,x_2)$) that $$R(y) = uP(y,x_2) + vxP^*(y,x_2) \equiv 0 \mod x^p.$$ The resultant $R(y)$ is a polynomial in $y$ with coefficients in $\overline{\mathbb{F}}_p[x]$; as the total degree of each of $P(y,x_2)$ and $xP^*(y,x_2)$ is at most $t$, it follows by definition of the resultant that \begin{equation*} \deg_y R(y) \leq \deg_y P(y,x_2) \cdot \deg_{x_2} xP^*(y,x_2) + \deg_y xP^*(y,x_2) \cdot \deg_{x_2} P(y,x_2) \leq 2t^2. \end{equation*} Similarly, the degree in $x$ of each coefficient in $R(y)$ is bounded from above by $2(s+1)t$. Therefore if $n$ is chosen with  $\max\{2t^2,2(s+1)t\} < n$ and $n^3 < p$, the conditions of Lemma 7 will be satisfied. It follows that $R(y) \equiv 0 \;\text{mod}\;x^p$.  Furthermore, as the degree in $x$ of each of the coefficients of $R(y)$ is less than $p$, it follows that $R(y)$ vanishes identically as a polynomial in $y$. 

We may assume that $t > 0$. After division by common factors of the coefficients of $P(x_1,x_2)$, we may assume that $P(x_1,x_2)$ is primitive in the sense of Gauss' lemma; as this only reduces the degree in $x$ of the coefficients, this does not interfere with the bounds required by the previous part of this proof. If $P(x_1,x_2)$ were reducible modulo $x^p$ in $\overline{\mathbb{F}}_p[x][x_1,x_2]$, then we could write \begin{equation}\label{eq24} P(x_1,x_2) \equiv f(x_1,x_2)g(x_1,x_2) \mod x^p,\end{equation} Any term in the product $f(x_1,x_2)g(x_1,x_2)$ containing a power of $x$ at least $p$ would originate from the product of powers $x^a$ and $x^b$ in $f$ and $g$, respectively, where one of $a$ or $b$ is at least $\frac{p}{2}$. Suppose without loss of generality that $a \geq \frac{p}{2}$. By setting the coefficient of $x^a$ in $f$ equal to zero, we again obtain \eqref{eq24}, as the degree $s$ of $P(x_1,x_2)$ in $x$ satisfies $s < \frac{n}{2t} < \frac{p}{2}$, and thus the term $x^a$ did not contribute to any nonzero term in $P(x_1,x_2)$. Therefore $P(x_1,x_2)$ would be reducible in $\overline{\mathbb{F}}_p[x][x_1,x_2]$ and by Gauss' lemma would thus be reducible in $\overline{\mathbb{F}}_p(x)[x_1,x_2]$.

Let us assume then that $P(x_1,x_2)$ is an irreducible polynomial. By the previous argument, it follows that $P(x_1,x_2)$ is irreducible modulo $x^p$. As $R(y)$ vanishes identically, it follows that $P$ and $xP^*$ are not coprime as polynomials in $x_2$. By irreducibility of $P$, it follows that $$xP^*(x_1,x_2) = T(x_1,x_2)P(x_1,x_2)$$ for some polynomial $T(x_1,x_2)$ with coefficients in $\overline{\mathbb{F}}_p(x)$. As the elements of $xP^*$ of a particular total degree in $x_1$ and $x_2$ derive from precisely those of the same total degree in $P$, it follows that $T$ has zero total degree in $x_1$ and $x_2$, and that $T=T(x)$ is an element of the rational function field $\overline{\mathbb{F}}_p(x)$. Let $H$ be the sum of terms in $P$ of highest total degree in $x_1$ and $x_2$. It follows that $$H^*(x_1,x_2) = \frac{T(x)}{x} H(x_1,x_2).$$ As in Siegel's argument \cite{Siegel}, it follows that there exists a nonzero solution $w$ to the differential equation \begin{equation} \label{eq25} w'' + \frac{1}{x} w' + w = 0 \end{equation} with $H\left(1,\frac{w'}{w}\right)= 0$. In particular, $u =\frac{w'}{w}$ is algebraic over $\overline{\mathbb{F}}_p(x)$ of degree at most $t < p$. Thus the extension $\overline{\mathbb{F}}_p(x)(u)|\overline{\mathbb{F}}_p(x)$ is tamely ramified. 

The function $u$ satisfies the Riccati differential equation \begin{equation} \label{eq26} u' + u^2 + \frac{1}{x}u = - 1. \end{equation} As ramification is tame in $\overline{\mathbb{F}}_p(x)(u)|\overline{\mathbb{F}}_p(x)$, we may write the Puiseux series for $u$ at a branch of infinity; as in Siegel's argument, we obtain by \eqref{eq26} that \begin{equation} \label{eq27} u = \pm 1 - \frac{1}{2x} + \cdots. \end{equation} Thus any branch of $u$ at infinity is regular and unramified. By \eqref{eq25}, the function $w$ is regular at all points $x \neq 0, \infty$, and thus branches of $u$ may only occur at zero or infinity. As infinity is not a branch point of $u$, it follows that zero is also not a branch point of $u$. Therefore $u$ is an element of $\overline{\mathbb{F}}_p(x)$. Similarly, expanding $u$ as a Laurent series at $x=0$ yields by \eqref{eq26} that $u$ is regular at zero. As $u$ is a rational function of $x$, the function $w$ has finitely many zeros, say $x_1,...,x_h$, and we may write \begin{equation} \label{eq28} u = \pm 1 + \sum_{k=1}^h \frac{1}{x - x_k}.\end{equation} By \eqref{eq27} and \eqref{eq28}, it follows that $h = -\frac{1}{2}$, a contradiction.

Our conditions on $A$, $B$, $C$, $D$, and $n$ in analogy to Theorems 1 and 2 may thus be written as \begin{enumerate} \item $n^3 < p$, \item $\max\{2(A+1)C,2C^2\} < n$, and \item $D(A+B+2D) < AB\cdot \frac{(C+2)}{2}.$ \end{enumerate} We set $n = \lfloor \frac{1}{2}p^{\frac{1}{3}} \rfloor$, which satisfies (1). For (2), we set $A = \lfloor \frac{1}{5}p^{\frac{2}{9}}\rfloor$ and $C = \lfloor  p^{\frac{1}{9}}\rfloor$. For (3), we set $D = \lfloor p^{\frac{2}{9}}\rfloor$ and $B = \lfloor 23p^{\frac{1}{9}}\rfloor$. Therefore $$N_{p}(y) \ll(A + pB + (p-1)C)/D \ll p^{\frac{8}{9}} .$$
 \end{proof}



\vspace{10pt}
{\small
\vspace{30pt}
\noindent AMIT GHOSH, Department of Mathematics, Oklahoma State University, Stillwater, OK 74078, USA \hfill {\itshape E-mail address}: ghosh@math.okstate.edu 
\vspace{12pt}

\noindent KENNETH WARD, Department of Mathematics, Oklahoma State University, Stillwater, OK 74078, USA \hfill {\itshape E-mail address}: kward@math.okstate.edu}

\begin{thebibliography}{9}
\footnotesize
\setlength{\parskip=0.0pt}
\setlength{\lineskip=0.0pt} 

\bibitem[Be08]{Beukers} F. Beukers, {\itshape E-functions and G-functions}, eprint available at  swc.math.arizona.edu/aws/2008/\\08BeukersNotesDraft.pdf .

\bibitem[Bo73]{Bo} E. Bombieri, {\itshape Counting Points on Curves over Finite Fields}, Seminaire Bourbaki, Springer, Berlin {\bf 430}, (1973) 234-241. 

\bibitem[Ch51]{Ch} C. Chevalley, {\itshape Introduction to the Theory of Algebraic Functions of One Variable}, Amer. Math. Soc. (1951).

\bibitem[Dw62]{Dw62} B. Dwork, {\itshape p-Adic cycles}, Pub. Math. I.H.E.S. vol. 37 (1969), 27-115.

\bibitem[He96]{He1} D. R. Heath-Brown, {\itshape An estimate for Heilbronn's exponential sum}, Analytic number theory: Proceedings of a conference in honor of Heini Halberstam, Birkh\"{a}user, Boston (1996) 451-463.

\bibitem[Ig58]{Ig58} J. Igusa {\itshape Class number of a definite quaternion with prime discriminant}, Proc. Nat. Acad. Sci. {\bf 44} (1958) 312-314.

\bibitem[Le87]{Le87} A. K. Lenstra Jr., {\itshape Factoring integers with elliptic curves}, Annals of Math. {\bf 126} No. 3 (1987) 649-673.

\bibitem[Mi92]{Mi} D. A. Mit'kin, {\itshape Stepanov method of the estimation of the number of roots of some equations}, Mat. Zametki {\bf 51} (1992), 52-58.



\bibitem[Si29]{Siegel} C. L. Siegel, {\itshape Uber einige Anwendungen diophantischer Approximationen}, Abhandlungen der K\"{o}niglichen Preu{\ss}ischen Akademie der Wissenschaften, Berlin (1929), 1-70.

\bibitem[St69]{St} S. A. Stepanov, {\itshape On the number of points of a hyperelliptic curve over a finite prime field}, Izv. Akad. Nauk SSSR {\bf 33} (1969) 1103-1111.

\end{thebibliography}
\end{document}